\documentclass[onefignum,onetabnum,final]{siamart250106}

\usepackage{geometry}
\usepackage{lipsum}
\usepackage{amsfonts}
\usepackage{graphicx}
\usepackage{epstopdf}
\ifpdf
  \DeclareGraphicsExtensions{.eps,.pdf,.png,.jpg}
\else
  \DeclareGraphicsExtensions{.eps}
\fi
\usepackage{lmodern}
\usepackage{amssymb,amsmath}
\usepackage{amsfonts}
\usepackage{algorithm}
\usepackage{algorithmicx}
\usepackage{algpseudocode}

\usepackage[nameinlink]{cleveref}
\crefname{line}{line}{lines}
\Crefname{line}{Line}{Lines}

\newcommand{\sign}{\mathrm{sign}}

\usepackage{calrsfs}
\DeclareMathAlphabet{\pazocal}{OMS}{zplm}{m}{n}
\usepackage{booktabs}
\usepackage{tikz}
\usepackage{pgf}  
\usepackage{subcaption} 

\newsiamremark{remark}{Remark}
\newsiamremark{hypothesis}{Hypothesis}
\crefname{hypothesis}{Hypothesis}{Hypotheses}
\newsiamthm{claim}{Claim}

\headers{Geometrical Adaptive Cross Approximation}{V.A. Yastrebov \& C. Noûs}

\title{Adaptive Cross Approximation with a \\Geometrical Pivot Choice: ACA-GP Method\thanks{Submitted to the editors \today}}

\author{Vladislav A. Yastrebov\thanks{Centre des matériaux, Mines Paris, PSL University, 91003 Evry, France
  (\email{vladislav.yastrebov@minesparis.psl.eu}, \href{https://yastrebov.fr}{www.yastrebov.fr}).}
\and Camille Noûs\thanks{Cogitamus Laboratory, Paris, France
  (\email{camille.nous@cogitamus.fr}, \href{https://www.cogitamus.fr/indexen.html}{www.cogitamus.fr}).}
}

\usepackage{amsopn}

\newcommand{\parencite}[1]{\cite{#1}}
\newcommand{\textcite}[1]{\cite{#1}}

\ifpdf
\hypersetup{
  pdftitle={Geometrical Adaptive Cross Approximation},
  pdfauthor={V.A. Yastrebov and C. Noûs}
}
\fi

\begin{document}

\maketitle

\begin{abstract}
    The Adaptive Cross Approximation (ACA) method is widely used to approximate admissible blocks of hierarchical matrices, or $\mathcal H$-matrices, from discretized operators in the boundary integral method. These matrices are fully populated, making their storage and manipulation resource-intensive. ACA constructs a low-rank approximation by evaluating only a few rows and columns of the original operator, significantly reducing computational costs. A key aspect of ACA's effectiveness is the selection of pivots, which are entries common to the evaluated row and column of the original matrix. This paper proposes combining the classical, purely algebraic ACA method with a geometrical pivot selection based on the central subsets and extreme property subsets. The method is named ACA-GP, GP stands for Geometrical Pivots. The superiority of the ACA-GP compared to the classical ACA is demonstrated using a classical Green operator for two clouds of interacting points.
\end{abstract}



\vspace{1em}
\noindent\textbf{Supplementary material:} The source code of the ACA and ACA-GP methods is available at \href{https://github.com/vyastreb/ACA}{github.com/vyastreb/ACA}, and the instance used for all studies presented in the paper is available at \href{https://github.com/vyastreb/ACA/releases/tag/v0.1.0}{github.com/vyastreb/ACA/releases/tag/v0.1.0} and \href{https://archive.softwareheritage.org/swh:1:dir:a9c10792ea93a27f342c938399f6cf4bf9916b98;origin=https://github.com/vyastreb/ACA;visit=swh:1:snp:c72b4194e2dcfb580ecdbf6a3e3caf4de37e3669;anchor=swh:1:rev:f03501b1dd63ba0658d43c435a1c2d122397c40f}{archive.softwareheritage.org}.
Supplementary material with all the data presented in the paper and all additional results is available at \href{https://doi.org/10.5281/zenodo.14809517}{zenodo.org/14809517}.
\normalsize
  
\normalsize
  
\tableofcontents

\section{Introduction}\label{introduction}

Low-rank approximation of matrices plays a crucial role in the numerical analysis of discretized differential or integral operators for boundary value problems in mathematical physics and various other domains. In the context of matrices arising from Boundary Integral (BIM) or Boundary Element (BEM) methods, they are often huge and dense, posing significant challenges in terms of storage and computational efficiency. Low-rank approximation techniques enable the reduction of these matrices to more manageable forms without sacrificing accuracy and improving efficiency of operations: various decompositions, matrix-vector product, inversion, etc. This simplification not only enhances computational performance but also makes it feasible to tackle with enhanced accuracy complex problems in scientific computing, engineering, and data science. Beyond mathematical physics, low-rank approximations are invaluable in fields such as machine learning, signal processing, and bioinformatics, where they facilitate efficient data compression, noise reduction, and the extraction of meaningful patterns from large datasets.

The method of Adaptive Cross Approximation (ACA) introduced in~\parencite{bebendorf2000approximation,bebendorf2003adaptive,bebendorf2006accelerating} has been successfully used to approximate admissible blocks of hierarchical matrices or $\mathcal H$-matrices~\parencite{bebendorf2008hierarchical} of discretized operators arising from the boundary integral method. Since the matrices associated with such operators are fully populated, their storage and manipulation are memory- and computationally intensive. The ACA permits not only to construct a low-rank approximation of the admissible blocks, but also it allows to control the approximation by evaluation only a few rows and columns of the original operator. This shortcut is especially beneficial when the computation of such entries is computationally expensive. The crucial point in an effective usage of the ACA is the choice of the so-called pivots which correspond to the common entry of the evaluated row and the column of the original matrix. In this article we suggest combining a purely algebraic classical ACA method with a geometrical choice of the pivots. This choice is based on the distance between the points in the cloud of points which represent the geometry of the problem. The proposed method is called ACA-GP (GP stands for Geometrical Pivots). We demonstrate the efficiency of the ACA-GP on a classical Green operator for elasticity.

The concept of $\mathcal H$-matrices was introduced in~\parencite{Hackbusch1999computing}, which is the most general class of hierarchical matrices with a strong admissibility. In $\mathcal H$-matrices, off-diagonal blocks can be low-rank and full-rank sub-matrices, a class of HODLR (hierarchical off-diagonal low-rank) matrices, as the name indicates, have only low-rank matrices at the off-diagonal blocks. So HODLR has a simpler structure.
The inverse of a matrix with low-rank blocks can be based on $\mathcal H$-matrix algebra~\parencite{bebendorf2008hierarchical,Grasedyck2003computing} or, as it was done in~\parencite{aminfar2016fast}, on the  Woodbury matrix identity
\[
(A + UCV)^{-1} = A^{-1} - A^{-1}U(C^{-1} + VA^{-1}U)^{-1}VA^{-1}
\]
which in the context of ACA decomposition has $C=I_{k\times k}$ and $C^{-1} = I$ and $UCV = UV^\top$.
HSS (Hierarchical Semiseparable) matrices~\parencite{chandrasekaran2008some,xia2010superfast} present a subclass of HODLR matrices and allow $\mathcal O(n)$-complexity solvers, an example can be found in~\parencite{Martinsson2005jcp}.
Another type of matrix structure is a Block Low-Rank (BLR) matrix format, a flat non-hierarchical block matrix structure~\parencite{amestoy2015improving}.
We can distinguish between a \emph{weak} and \emph{strong admissibility} of non-overlapping domains as was introduced in~\parencite{ambikasaran2013fast}.

Exist numerous techniques to construct low-rank approximations of admissible blocks $K \in\mathbb R^{m\times n}$, $n<m$. They could be subdivided in two main classes: (i) a very general algebraic skeleton-based or pseudo-skeleton approximations~\parencite{goreinov1997theory}, where particular rows and columns of the original matrix are selected, (ii) by direct decompositions of the matrix into more adapted for low-rank approximation forms such as SVD, LU and QR, and (iii) kernel approximations, where the underlying operator of two variables is approximated as a product of two separate operators acting on one variable. The first class is the most efficient in terms of computational resources, the second is the most general and accurate, the two can be applied to any matrix, the third class is more specific and is applicable when the matrix is defined by a smooth kernel operator. A survey could be found in~\parencite{grasedyck2013survey}.
%
%
The most popular methods are summarized below:
\begin{itemize}
  \item \textbf{Interpolative Decomposition (ID)} based on matrix skeletonization~\parencite{cheng2005compression}, complexity $\mathcal O(mn\log k+k^2n)$ or as described in~\parencite{woolfe2007fast} as at most $\mathcal O(kmn\log(n))$ but typically $\mathcal O(kmn)$.
  \item \textbf{Randomized algorithms}~\parencite{frieze2004fast,deshpande2006adaptive,halko2011finding} with a fast randomized algorithm~\parencite{woolfe2007fast,woolfe2008fast} with the complexity $\mathcal O(nm log(k) + nl^2)$ where $l$ is of order $k$ but $k<l<n$. The fast version is based on the ID and requires a knowledge of the matrix
  \item \textbf{Adaptive Cross Appoximation (ACA)} with full pivoting needs to scan the full matrix and Partial pivoting Adaptive Cross Approximation (ACA) and ACA+, seek for largest entries in a row or column at every approximation step; they work well for rather homogeneous matrix structure, but can fail when especial rows and columns are present in the matrix, complexity $\mathcal O(k(m+n))$.
  \item \textbf{Boundary Distance Low-Rank (BDLR)}~\parencite{aminfar2016fast}, relies on the underlying sparse matrix graph to choose the desired rows and columns, complexity $\mathcal O(nk)$ however it can be used in a special context and requires the original sparse matrix from the finite element solver.
  \item \textbf{CUR} decomposition~\parencite{mahoney2009cur} at the boundary between classical decompositions and pseudo-skeleton approximations, where the selected rows and columns exhivit a "high statistical leverage", but it is mainly used in the context of data compression and is not directly related to matrices obtained from BIM, which is the main focus of this work.
  \item \textbf{Rank revealing $QR$} algorithm with a pivoted Gram-Schmidt algorithm to obtain $r$ orthogonal basis vectors, used, for example, in~\parencite{kong2011adaptive}, complexity $\mathcal O(mnk)$
  \item \textbf{Rank revealing $LU$} decomposition, complexity $\mathcal O(mnk)$
  \item \textbf{Singular Value Decomposition} (SVD)~\cite[e.g.]{golub2013matrix}, complexity for a full SVD construction $\mathcal O(mn^2 + n^3)$ and for a truncated SVD $\mathcal O(kmn)$.
  \item \textbf{Chebyshev low-rank approximation} for smooth kernels is an expansion of the kernel as a series of Chebyshev polynomials, complexity $\mathcal O(nk)$.
  \item \textbf{Multipole expansion} is another type of expansion of smooth kernels as series of spherical harmonics, it is used in the Fast Multipole Method~\parencite{greengard1987fast}, complexity $\mathcal O(kn)$. 
\end{itemize}

In this short paper we develop a new method for a low-rank approximation of matrices arising from physical interaction between admissible (well separated) domains described by a smooth kernel operator. It is based on the psedo-skeleton approximation and is largely inspired by a purely algebraic partial pivot ACA method, but the choice of pivots is based on the combination of algebraic and geometrical properties of the problem. This paper is organized as follows. In \Cref{sec:methods} we present the ACA method and its proposed ACA-GP modification for point-wise interaction. In \Cref{sec:results} we present the results of the ACA-GP method applied to a Green operator for a set of separate clouds of points. 
The performance of the method is compared with the classical ACA and the most computationally extensive but the most accurate Singular Value Decomposition (SVD) ensuring the best low-rank approximation. 
Finally, in \Cref{sec:conclusions} we draw some conclusions and propose future research directions.

\section{Methods\label{sec:methods}}

\subsection{Low rank approximation of admissible blocks}

Let us consider admissible blocks of a hierarchical matrix~\parencite{bebendorf2008hierarchical} of a discretized operator arising from the 
Boundary Integral Method (BIM)~\parencite{bonnet1999boundary}. 
By admissible blocks we understand the blocks corresponding to subdomains $X,Y\subset\Omega$ which verify the following geometrical criterion~\parencite{bebendorf2008hierarchical}:
\begin{equation}
  \min\{\mathrm{diam}(X), \mathrm{diam}(Y)\} \le \eta \mathrm{dist}(X,Y),
  \label{eq:admissibility}
\end{equation}
where by $\mathrm{diam}(\bullet)$ we understand a certain measure of the set extension which can be precomputed with a reasonable (linear) complexity, for example as a the doubled distance from the center to the farthest element. Since the computation of the minimal distance between the sets has a high ($\mathcal O(nm)$) complexity, a "relaxed" distance condition shall be used.  This relaxed distance will be denoted with a prime $\mathrm{dist'}(X,Y)$, meaning that this is an upper bound of the real Euclidean distance. This approximate distance can be also computed with a reasonable (linear) complexity, for example, as a distance between centers of subdomains $\bar x,\bar y$ of $X$ and $Y$ minus the estimated half-sum of diameters, i.e.
\begin{equation*}
\begin{split}
  \mathrm{dist(X,Y)} &\le \mathrm{dist'}(X,Y) = 
  \|\bar x - \bar y\| - (\mathrm{diam}(X) - \mathrm{diam}(Y))/2 \le \\ &\le \|\bar x - \bar y\| - \min\{\mathrm{diam}(X), \mathrm{diam}(Y)\}.
\end{split}  
\end{equation*}
For a set of points, the center can be computed with a linear complexity as a barycenter:
\begin{equation}
  \bar x = \frac{1}{|X|} \sum_{i\in X} x_i, \quad \bar y = \frac{1}{|Y|} \sum_{j\in Y} y_j,
\end{equation}
where $|X|,|Y|$ denote the number of points in the sets $X,Y$, respectively.
Then the admissibility criterion \Cref{eq:admissibility} reduces to
\begin{equation}
  \min\{\mathrm{diam}(X), \mathrm{diam}(Y)\} \le \alpha \|\bar x -\bar y\|, \quad \alpha = \frac{\eta}{1+\eta}.
\end{equation}

Now let us assume that the components $a_{ij}$ of the admissible block of the discrete BIM operator $A \in \mathbb R^{n\times m}$ can be computed as
$$
 a_{ij} = \int\limits_{X_i} \int\limits_{Y_j} \kappa(x,y) \, dx \, dy, \quad i=1,\ldots,n, \quad j=1,\ldots,m
$$
where $\kappa(x,y)$ is the kernel operator. Since the subdomains are separated, there is no singular function evaluation. The objective is to evaluate a low-rank approximation $A'_k$ of the operator $A$ such that the rank $k$ of the approximation ensures a controllable error $\|A-A'_k\|_F$ and that $k$ preferably remains much smaller than the size of the original matrix $k \ll \min\{m,n\}$. 

\subsection{Adaptive Cross Approximation (ACA)}

The ACA method~\parencite{bebendorf2003adaptive} is based on the general idea of a pseudo-skeleton representation of the matrix $A$~\parencite{goreinov1997theory}.
The approximate low-rank matrix $A'_k$ of rank $k$ is represented as an outer product of submatrices $U$ and $V$:
\begin{equation}
  A'_k = \sum_{i=1}^k u_i v_i^T = U V^T, \qquad U \in \mathbb R^{n\times k}, V \in \mathbb R^{m\times k}
\end{equation}
which can be adaptively constructed until a required accuracy is reached:
\begin{equation}
  \forall \varepsilon > 0, \exists k \quad \text{s.t.} \quad\|A - A'_k\|_F \le \varepsilon \|A\|_F,
\end{equation}
where $\varepsilon$ is the required accuracy $\epsilon$ and $\|\bullet\|_F$ is the Frobenius norm. Of course, since the matrix $A$ is never computed fully nor the low-rank matrix $A'$  is constructed (only its decomposed representation), all these norms are replaced by adapted expressions and will be presented in the following section.
The data compression ratio is defined as
\[
  \frac{|U|+|V|}{n \times m} = \frac{k(n+m)}{n \times m} \xrightarrow[]{m=n} \frac{2k}{n}
\]
which clearly demonstrates that this compression makes sense only if $k\!<\!\min\{n,m\}/2$. 

The classical ACA method is based on the algorithm shown in \Cref{alg:aca}. The only liberty in the algorithm is the choice of the pivot row $i_k$ (line 3)
from which the pivot column $j_k$ is selected such that it has the maximal entry of the residual column (line 10).
The choice of the pivots is crucial for the convergence of the ACA method and classically the row indices $i_k$ are selected such that the Vandermonde matrix corresponding to the system in which the approximation error is to be estimated is non-singular, see~\parencite{bebendorf2008hierarchical,bauer2021block}. In~\parencite{grasedyck2005adaptive} a slightly different strategy is applied and called ACA+: starting with a random initial row and subsequent maximal component choice for the column and the maximal entry of this column for the row. The ACA+ method is claimed to be more robust than the classical ACA but as demonstrated in~\parencite{bebendorf2008hierarchical}, in some situations it could be suboptimal.

Let us now demonstrate the ACA type matrix approximation in a more explicit form assuming that components $a_{ij}$ are computed as $a_{ij} = 1/|x_i - y_j|$. 
The first rank $k=1$ is constructed as
$$
{A'_1}_{ij} = \frac{a_{i_1,j} a_{i,j_1}}{a_{i_1,j_1}} = \frac{|x_{i_1}-y_{j_1}|}{|x_i - y_{j_1}|\,|y_j - x_{i_1}|}, \quad i=1,\ldots,n, \quad j=1,\ldots,m
$$
where ${i_1,j_1}$ is the first pivot.
The associated residual between the original matrix and this first rank approximation is given by
\begin{equation}
  R^{1}_{ij} = A_{ij} - {A'_1}_{ij} = \frac{1}{|x_i - y_j|} - \frac{|x_{i_1}-y_{j_1}|}{|x_i - y_{j_1}|\,|y_j - x_{i_1}|}.
\label{eq:R1}
\end{equation}
We can represent points in the clouds with respect to the first pivot as $x_i = x_{i_1} + \Delta^1 x_i$ and $y_j = y_{j_1} + \Delta^1 y_j$, where $\Delta^1 x_{i_1} = 0$ and $\Delta^1 y_{j_1} = 0$.
If we introduce the following notation for the vector connecting two pivot points $d_{11} = x_{i_1} - y_{j_1}$, then the residual becomes
\begin{equation}
R^{1}_{ij} = 
\frac{|d_{11} + \Delta^1 x_i||d_{11} - \Delta^1 y_j| - |d_{11} + \Delta^1 x_i - \Delta^1 y_j||d_{11}|}{|d_{11} + \Delta^1 x_i - \Delta^1 y_j||d_{11} + \Delta^1 x_i||d_{11} - \Delta^1 y_j|} 
\label{eq:R1d}
\end{equation}
It is easy to see that for $i=i_1$ and for $j=j_1$ the numerator is equal to zero, i.e. the residual column and line corresponding to the pivot are zero.
The second rank $k=2$ approximation is constructed as
$$
{A'_2}_{ij} = {A'_1}_{ij} + \frac{R^{1}_{i,j_2} R^{1}_{i_2,j}}{R^{1}_{i_2,j_2}}
$$
where ${i_2,j_2}$ is the second pivot. The residual between the original matrix and this second rank approximation is given by
$$
R^2_{ij} = A_{ij} - {A'_2}_{ij} = R^1_{ij} - \frac{R^{1}_{i,j_2} R^{1}_{i_2,j}}{R^{1}_{i_2,j_2}} = \frac{R^{1}_{i,j} R^{1}_{i_2,j_2} - R^{1}_{i,j_2} R^{1}_{i_2,j}}{R^{1}_{i_2,j_2}}.
$$
With this form it is easy to see that for $i=\{i_1,i_2\}$ and $j=\{j_1,j_2\}$ the residual is zero.
Recursively, the $(k+1)$ rank approximation is constructed as
$$
R^{k+1}_{ij} = A_{ij} - {A'_{k+1}}_{ij} = R^k_{ij} - \frac{R^k_{i,j_{k+1}} R^k_{i_{k+1},j}}{R^k_{i_{k+1},j_{k+1}}} = \frac{R^k_{i,j} R^k_{i_{k+1},j_{k+1}} - R^k_{i,j_{k+1}} R^k_{i_{k+1},j}}{R^k_{i_{k+1},j_{k+1}}}.
$$

\begin{algorithm}
  \caption{Adaptive Cross Approximation (ACA)\label{alg:aca}}
  \begin{algorithmic}[1]
    \State \textbf{Let} $k \gets 1$; $I,J \gets \emptyset$; $\epsilon > 0$
    \Repeat
        \State Find pivot row $i_k$ by some rule\label{alg:aca:pivot}
        \State Evaluate column $\tilde{v}_k \gets A_{i_k,s}$
        \For{$l = 1, \ldots, k - 1$} \Comment Subtract previous columns
            \State $\tilde{v}_k \gets \tilde{v}_k - (u_l)_{i_k} v_l$
        \EndFor
        \State $I \gets I \cup \{i_k\}$
        \If{$\tilde{v}_k$ does not vanish}
            \State $j_k \gets \arg \max_{j \in S} |\tilde{v}_k|_j$\label{alg:aca:pivot2}
            \State $J \gets J \cup \{j_k\}$
            \State Pivot $p_k \gets \tilde{v_k}_{j_k}$
            \State Evaluate row $\tilde{u}_k \gets A_{t,j_k}$
            \For{$l = 1, \ldots, k - 1$} \Comment Subtract previous rows
                \State $\tilde{u}_k \gets \tilde{u}_k - (v_l)_{j_k} u_l$
            \EndFor
            \State Evaluate residual norm $\|R_k\|_F$ and matrix norm $\|A'_k\|_F$
            \State Renormalize $u_k = \mathrm{sign}(p_k)\tilde{u}_k / \sqrt{p_k}$, $v_k = \tilde{v}_k / \sqrt{p_k}$
            \State Update matrices $U \gets [U|u_k]$, $V \gets [V|v_k]$
            \State $k \gets k + 1$
        \EndIf
    \Until{$\|R_k\|_F \leq \epsilon \|A'_k\|_F$ or $|I| = n$ or $|J| = m$}
  \end{algorithmic}
\end{algorithm}

\subsection{Geometrical insights in optimal pivot choice}

The exclusively algebraic search for pivots used in ACA, based on maximizing pivot values, has the merit of being geometry-agnostic and easy to implement.
However, it is suboptimal, as can be easily demonstrated through a relatively simple test. 
Take two random separate clouds $X$ and $Y$ of size $|X|=N$, $|Y|=M$ and compute all possible choices for the first pivot $\{i_1,j_1\}$ such that $\{i_1,j_1\} \in N \times M = \{1,2,\dots,n\} \times \{1,2,\dots,m\}$ for every pivot we can construct vectors $\tilde u_1 = a_{i,j_1}$ and $\tilde v_1 = a_{i_1,j}$ and approximate the matrix by 
$$
A'_1(i_1,j_1) = \tilde u_1 \tilde v_1^\top / p_{i_1,j_1}, \quad p_{i_1,j_1} = a_{i_1,j_1}.
$$
Then we can select the optimal pivot $\{i_1^*,j_1^*\}$ as
$$
\{i_1^*,j_1^*\} = \arg\min\limits_{\{i_1,j_1\}\in N \times M} \left(\|A'_1(i_1,j_1) - A\|_F\right),
$$
by definition 
$$
\|A'_1(i_1^*,j_1^*) - A\|_F \le \|A'_1(i_1,j_1) - A\|_F, \forall i_1 \in N,j_1 \in M.
$$
In the ACA $i_1$ is selected randomly and $j_1 = \arg\max_{j\in M} |A_{i_1,j}|$. 
We can visualize the error in the approximation for the two clouds for a fixed row or a column, i.e. $\forall i, \{i,j_1\}$ for cloud $X$ and $\forall j, \{i_1,j\}$ for cloud $Y$. It makes sense to show the error for $i_1=i_1^*$ and $j_1=j_1^*$. Now we proceed to the error estimation of the next rank $k$ by assuming that on the previous rank we selected the optimal pivot $\{i_{k-1}^*,j_{k-1}^*\}$. The scaled rows and columns are stored in matrices $U = \{u_1, u_2, \dots, u_k\}$ and $V = \{v_1, v_2, \dots, v_k\}$, with 
\begin{align}
  u_l &= \sign(p(i_l,j_l))\tilde u_l / \sqrt{|p(i_l,j_l)|},\\
  v_l &= \tilde v_l / \sqrt{|p(i_l,j_l)|}
\end{align}
so that $A'_k = U V^\top$, where $\tilde u_l,\tilde v_l$ are the row and the column vectors of the residual matrix $R_{k-1}$ and $p(i_l,j_l)$ is the pivot value:
\begin{align}
  \tilde u_l &= R^{l-1}_{i,j_l} = a_{i,j_l} - a'^{l-1}_{i,j_{l-1}},\\
  \tilde v_l &= R^{l-1}_{i_l,j} = a_{i_l,j} - a'^{l-1}_{i_{l-1},j},\\
  p(i_l,j_l) &= \tilde u_l(i_l) = \tilde v_l(j_l).
\end{align}
 The relative error for the pivot choice $\{i_l,j_l\}$ can be computed as
$$
E_{k}(i_{k},j_{k}) = \frac{\|A'_{k}(i_{k},j_{k}) - A\|_F}{\|A\|_F}.
$$
Again, we can find the optimal choice for the pivot
$$
\{i_{k}^*,j_{k}^*\} = \arg\min\limits_{\{i_{k},j_{k}\}\in \{N\setminus I\} \times \{M \setminus J\}} E_{k}(i_{k},j_{k}),
$$
where $I=\{i_1,i_2,\dots,i_{k-1}\}$, $J=\{j_1,j_2,\dots,j_{k-1}\}$ 
 contains rows and columns of previously selected pivots, thus $N \setminus I$ and $M \setminus J$ are sets of the remaining rows and columns. Again we can visualize the error in two clouds for a fixed row or a column, i.e. $\forall i, \{i,j_{k}^*\}$ for cloud $X$ and $\forall j, \{i_{k}^*,j\}$ for cloud $Y$. 
Let's call this search strategy \textit{genetic} as for every rank we select the best possible approximation, which is of course, unimaginable in practice as this requires construction of $n\times m$ full matrices $n\times m$ times to find the most optimal one, so, the complexity would be $\mathcal O(kn^2m^2)$ which is too prohibitive in practice. But here our objective is first to demonstrate the sub-optimality of the classical ACA pivot selection, and second 
to reveal the geometrical structure of the optimal pivot location. By pivot location we mean the geometrical coordinates of the optimal row $i_k$ (point $x_{i_k}$) and column $j_k$ (point $y_{j_k}$) in the clouds $X$ and $Y$, respectively. 

\begin{figure}[h!]
  \centering\includegraphics[width=1\textwidth]{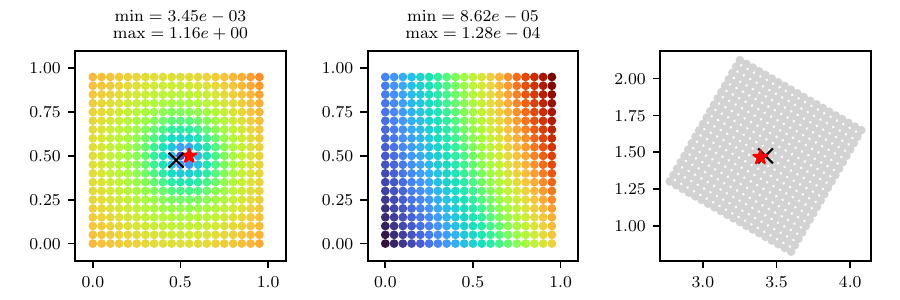}
  \caption{\label{fig:rank1geometry}Illustration of the optimal choice of the first pivot $\{i_1^*,j_1^*\}$ for the rank $k=1$ for two clouds of points ($x_{i_1^*},y_{j_1^*}$ are marked with stars). The optimal pivot points are located near the centers of the clouds. Left column: relative error above the SVD approximation: $\tilde E_1 = (E_1 - E_1^{\text{\tiny{SVD}}})/E_1^{\text{\tiny{SVD}}}$. Central column: interaction matrix's column $a_{i_1,j}$ colored according to their value. Right column: the location of the optimal column $i_1$ in the cloud $X$. Black crosses correspond to the barycenters of the clouds.}
\end{figure}

To analyse the geometrical signature of the optimal pivot choice we analyse two clouds of points $X$ and $Y$ located at some distance and oriented in a random fashion. Random and structural point distributions within the clouds are studied based on the full genetic data generated.
To make the error estimation more meaningful, we subtract the relative error of the SVD approximation for the corresponding rank $k$ and normalize by the relative SVD error:
\begin{equation}
\tilde E_k = \frac{E_k - E_k^{\text{\tiny{SVD}}}}{E_k^{\text{\tiny{SVD}}}} = \frac{\|A'_{k} - A\|_F}{\|A'^{\text{\tiny{SVD}}}_{k} - A\|_F} - 1\ge 0,
\label{eq:tildeE}
\end{equation}
since, by definition the SVD is the optimal low-rank approximation, we can be sure that the error remains positive.

Based on the data analysis of the genetic ACA, the \emph{first observation} we make is that the optimal choice of the first pivot is located near centers of the clouds $X$ and $Y$, which is quite obvious. It can be readily seen in \Cref{fig:rank1geometry} that the points near the center possess the minimal error.

\begin{figure}[ht!]
  \centering\includegraphics[width=0.9\textwidth]{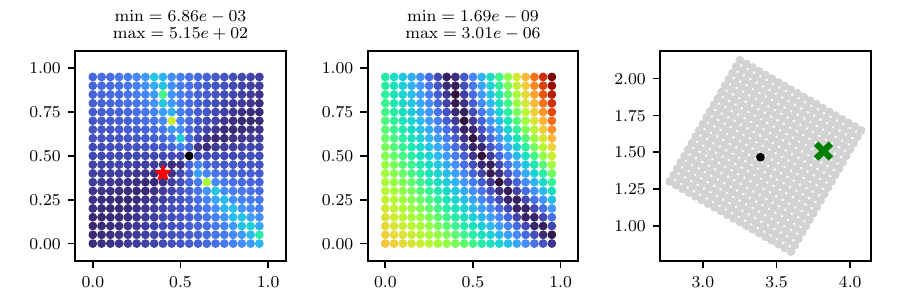}
  \centering\includegraphics[width=0.9\textwidth]{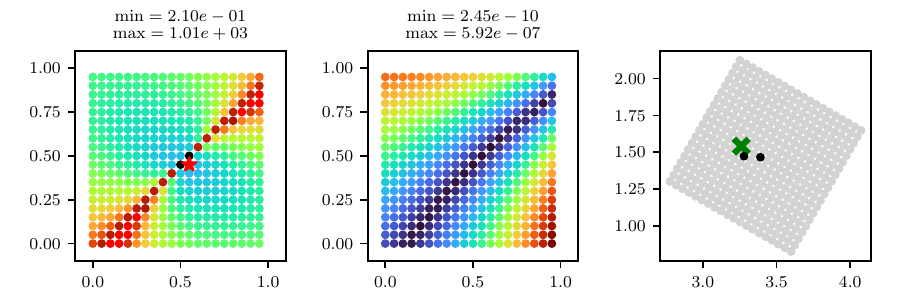}
  \centering\includegraphics[width=0.9\textwidth]{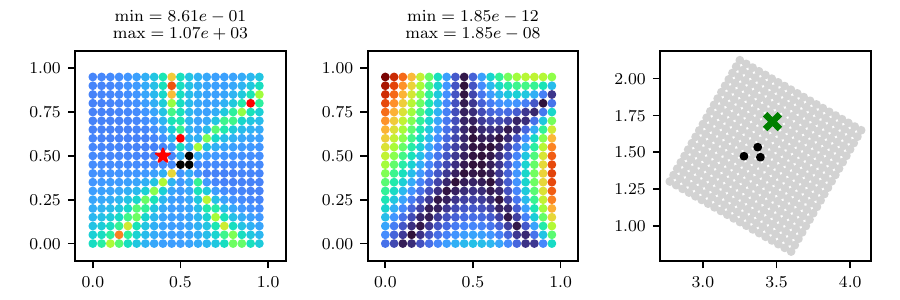}
  \centering\includegraphics[width=0.9\textwidth]{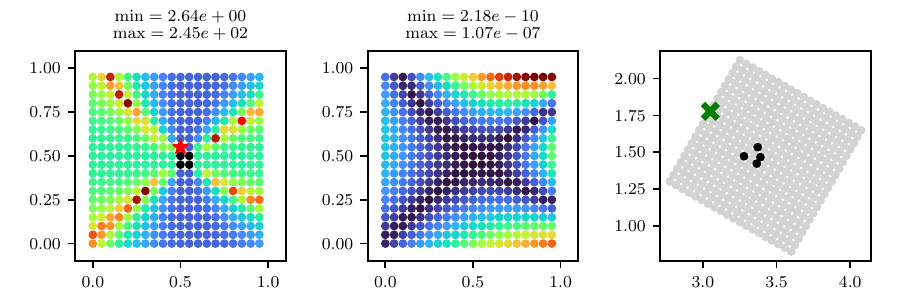}
  \caption{Illustration of the error structure for ranks $k=\{2,3,4,5\}$ starting from the top and a randomly selected $\{i_2,i_3,i_4,i_5\}$ (note that all pivots $i_l^*,j_l^*$ with $l<k$ were assumed to be selected in an optimal way and are shown with black circles) in cloud $X$ shown with a cross on the right column; the left column represents the relative difference in genetic ACA's error (relative error of the Frobenius norm of the matrix approximation) with respect to SVD $\tilde E_k = (E_g - E_{\text{SVD}})/E_{\text{SVD}}$; the central column shows the absolute value of the residual column corresponding to selected $i_k$; the clouds have $n=m=400$ points.}
  \label{fig:rank234}
\end{figure}

The \emph{second observation} is that starting from rank $k=2$, minimal absolute values of the residual matrix's column corresponds to the maximal approximation errors. By residual matrix, as previously, we understand $R_{k-1} = A-A'_{k-1} = a_{ij} - \sum_{i=1}^{k-1} u_i v_i^\top$. So for every $i_k^*$ we can evaluate all values $R_{k-1}(i_k^*,j)$; in the classical ACA we select the maximal value $j_k^* = \arg_j\max|R_{k-1}(i_k^*,j)|$, but as can be seen for $k=1$ in \Cref{fig:rank1geometry} (for which $R_0 = A$), and for higher ranks as illustrated in \Cref{fig:rank234}, this choice is not optimal. Nevertheless, what is clear is that the minimal absolute values of the residual column correspond to the maximal errors of the approximation $\tilde E_k$ and thus the ACA by choosing the higher values of the residual safely avoids too bad approximation. 

The \emph{third observation} is that the optimal points concentrate in the center of the cloud. It could be readily seen in \Cref{fig:rank234}: black circles correspond to the optimal pivots selected using genetic ACA. It contrasts with the classical ACA choice in which the maximal values of the residual correspond to points on the boundary of the cloud as can be easily seen in the central column of \Cref{fig:rank234}.

The \emph{fourth observation} is that the zones of extreme errors (minimal or maximal) lie on some curved lines which seem to pass through optimal pivots of the previous ranks. These curves, let's call them \emph{extreme curves}, could be approximated by some function $\mathcal E_k: f(x) = 0$ such that error minimizing and maximizing pivots lie on these extreme curves: $f(x_{i}) \approx 0 \approx f(y_{j})$. The sign $\approx$ reflects the discrete nature of points and the continuous nature of the curves. We can distinguish extreme curves of maximal $\mathcal E_k^{\max}$ and minimal $\mathcal E_k^{\min}$ errors. This observation is valid for $k > 1$. 
Regardless a seeming simplicity of these extreme curves and relatively simple origin algebraic equation for their construction, they cannot be easily determined for arbitrary domains. Therefore, the construction of the algorithm will combine such curves for ranks $k=2$ and $k=3$ with the third observation about the optimal points concentration in the center of the cloud.
Nevertheless, we will make an attempt to reveal the geometrical structure of these extreme curves for simple domains and formalize the algorithm for the optimal pivot selection in the following sections.

\subsection{Geometrical structure of the error distribution\label{sec:central_subsets}}

To take profit of the observation that the optimal points concentrate in the center of the cloud, we introduce two central subsets $I^c,J^c$ of the clouds such that
\begin{align}
  I^c &= \left\{ i \in \{1,\dots,n\} \left| \vphantom{\frac A A} \; |x_{i} - x_{i_1}| \le \varepsilon_r \mathrm{diam}(X)\right.\right\}\\
  J^c &= \left\{ j \in \{1,\dots,m\} \left| \vphantom{\frac A A} \; |y_{j} - y_{j_1}| \le \varepsilon_r \mathrm{diam}(Y)\right.\right\},
  \label{eq:central_subsets}
\end{align}
which contain indices of points close to the points forming the first pivot $x_{i_1}, y_{j_1}$. The parameter $\varepsilon_r$ is a relative distance to the diameter of the cloud. The only condition on the selection of the central subsets is that they must contain sufficient number of points for the relevant range of ranks, for example, if the maximal rank is bounded by $\max(\text{rank})$ then the number of points in the central subset must contain at least this number with a relatively small offset of the order of $\Delta \approx 5-10$: $|I_c| \ge \max(\text{rank}) + \Delta$. The algorithm for the central subset selection is presented in \Cref{alg:central_subset}. If the required number of points is not verified, the central ratio fraction is increased, for example, by $10$ \%.

The sensitivity of the approximation to this central fraction -- the only parameter of the algorithm -- will be evaluated in \Cref{sec:results}.
The central subsets $I^c,J^c$ will be used to construct \emph{ad hoc} rules for the selection of optimal points. 
As will be shown, a selection algorithm based on extreme curves can be constructed for more or less square domains for the ranks $k=2$ and $k=3$, for which the geometrical structure of error is rather trivial. 
For higher ranks $k\ge4$ or arbitrary domains a simpler procedure combining central subsets and algebraic considerations is adopted, which will be presented in \Cref{sec:algorithm_for_high_ranks}.

\begin{algorithm}
  \caption{Central subset selection\label{alg:central_subset}}
  \begin{algorithmic}[1]
    \State \textbf{Let} first pivot $\{i_1\}$ be selected
    \State \textbf{Let} $I^c \gets \emptyset$; $\varepsilon_r > 0$
    \While {$|I^c| < \max(\text{rank})+\Delta$}
    \For{$i = 1, \ldots, n$}
        \If{$|x_i - x_{i_1}| \le \varepsilon_r \mathrm{diam}(X)$}
            \State $I^c \gets I^c \cup \{i\}$
        \EndIf
    \EndFor
    \If{$|I^c| < \max(\text{rank})+\Delta$}
    \State $\varepsilon_r \gets 1.1 \varepsilon_r$
    \EndIf
    \EndWhile
  \end{algorithmic}
\end{algorithm}

\subsubsection{Rank $k=2$}

Let us start the study from the second low-rank approximation $k=2$, we assume that we selected the optimal pivot $\{i^*_1,j^*_1\}$ for the first rank. The objective is to select such point $x_{i^*_2}$ in cloud $X$ and such point $y_{j^*_2}$ in cloud $Y$ so that the newly approximated matrix $A'_2 = (u_1,u_2) (v_1,v_2)^\top$ minimizes the Frobenius norm $\{i_2^*,j^*_2\}=\arg\min_{i,j}\|A-A'_2\|_F$. 
If a point in cloud $X$ is selected in random fashion $i_2$, then there exist such $j'_2$ which minimizes the aforementioned norm for given $i_2$. At the same time, for every choice of $i_2$, one can visualize the error $\tilde E_2(i_2,j)$, see \Cref{eq:tildeE}, for all points of the cloud $Y$ as it was done before. 

\begin{figure}[ht!]
  \includegraphics[width=1\textwidth]{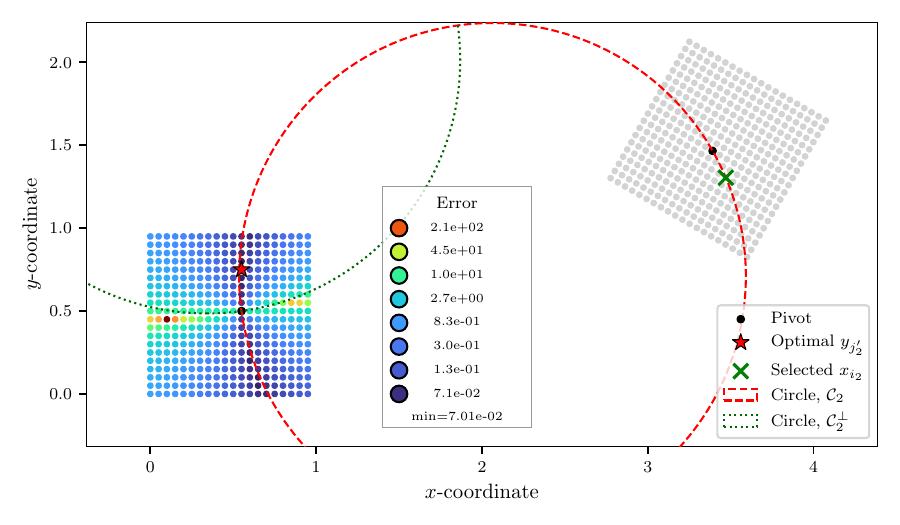}
  \includegraphics[width=1\textwidth]{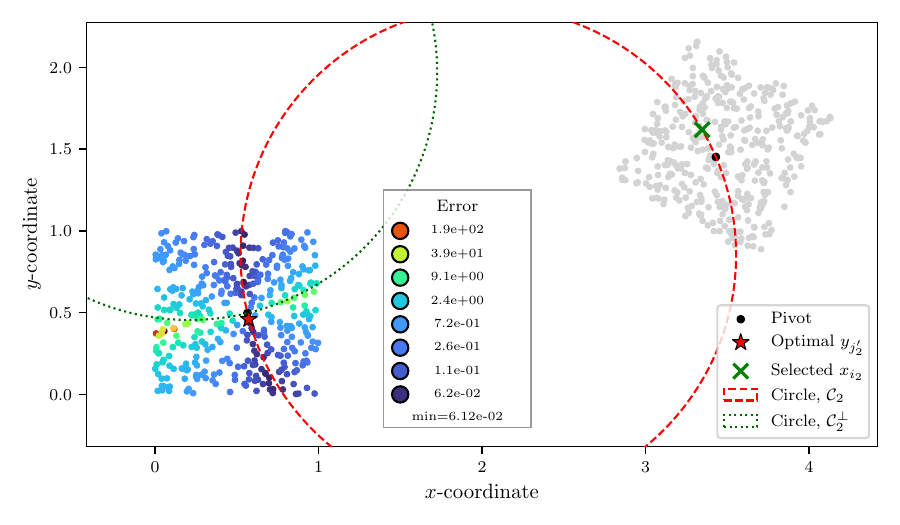}
  \caption{\label{fig:error_k2a}Illustration of the error structure for the rank $k=2$ for two clouds of structured (upper panel) and random (lower panel) points within square boxes. A randomly selected point $x_{i_2}$ in $X$ cloud is marked with a green cross and the corresponding optimal $y_{j'_2}$ is shown with a star, the corresponding circle $\mathcal C_2(i_2)$ is shown with a red dashed line, the conjugate circle $\mathcal C_2^{\perp}$ is shown with a green dotted line; black circles highlight the first optimal pivot $x(i_1^*)$ and $y(j_1^*)$.}
\end{figure}

\begin{figure}[ht!]
  \includegraphics[width=1\textwidth]{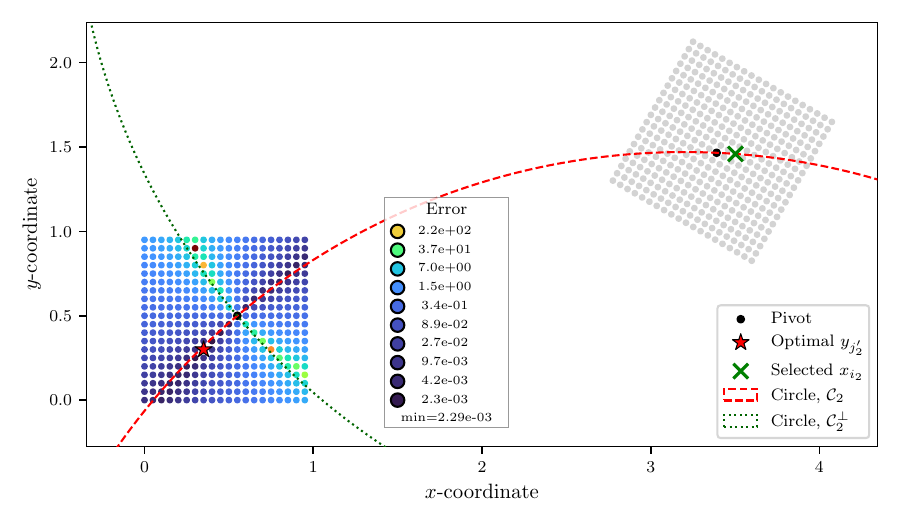}
  \includegraphics[width=1\textwidth]{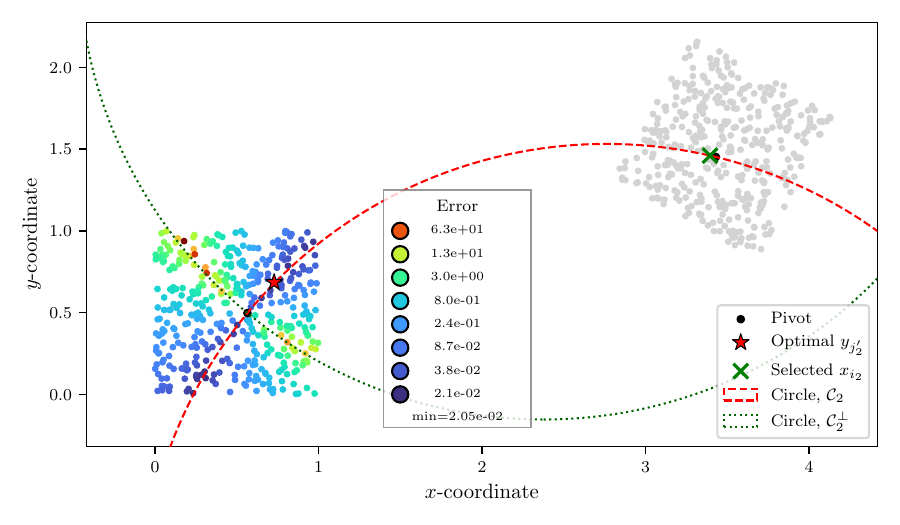}
  \caption{\label{fig:error_k2b}Illustration of the error structure for the rank $k=2$ for two clouds of structured (upper panel) and random (lower panel) points within square boxes. A randomly selected point $x_{i_2}$ in $X$ cloud is marked with a green cross and the corresponding optimal $y_{j'_2}$ is shown with a star, the corresponding circle $\mathcal C_2(i_2)$ is shown with a red dashed line, the conjugate circle $\mathcal C_2^{\perp}$ is shown with a green dotted line; black circles highlight the first optimal pivot $x(i_1^*)$ and $y(j_1^*)$.}
\end{figure}

On the top of the clouds colored according to relative error $\tilde E_k$, we make a simple geometrical construction: a circle $\mathcal C_2(i_2)$ passing through points $x_{i^*_1}$, $y_{j^*_1}$ and the new (tentative) point $x_{i_2}$. According to our experiments, for any choice of $i_2$, the optimal $y(j'_2)$ lies on the same circle with an eventual spatial shift related to the spatial distribution of points. Therefore, the extreme curve for $i_2$: $\mathcal E_2^{\min}(i_2)$ can be approximated by $\mathcal C_2(i_2)$. Let us construct a \emph{conjugate circle} passing through the point $y_{j^*_1}$: $\mathcal C_2^{\perp}(y_{j^*_1})$ (or for shorter notation $\mathcal C_2^{\perp}(j^*_1)$) in such a way that it has the same radius as $\mathcal C_2$, it intersects the point $y_{j^*_1}$ orthogonally to the circle $\mathcal C_2$, and the center $c_2^{\perp}$ verifies $(c_2^{\perp} - y_{j^*_1})\cdot(x_{i_1^*}-y_{j_1^*}) \ge 0$. Constructed in such a way, this conjugate circle represents the extreme curve of maximal errors $\mathcal E^{\max}_2(i_2) = \mathcal C^{\perp}_2$. Interestingly, regularly structured clouds of points possess the same geometrical structure of the error as randomly distributed clouds with a relatively high and uniform density. In \Cref{fig:error_k2a,fig:error_k2b} we show the error distribution for randomly selected $i_2$ along with the introduced geometrical constructions: the circle $\mathcal C_2(i_2)$ and the conjugate circle $\mathcal C_2^{\perp}(y_{j^*_1})$. Only two configurations are shown, more examples are provided in Supplementary material~\parencite{supplementary}.

Now, we need to develop a set of \emph{ad hoc} rules which would help us to select more or less optimal points. In order to keep these rules simple to implement, we propose an algorithm shown in \Cref{alg:rank2_selection}.
It aims to iteratively select an optimal pair of points from two clouds to keep small the approximation error for a rank-2 matrix. The process starts by randomly selecting a point from one of the clouds $x_{i_2}\in X$ within the central subset $i_2\in I^c\setminus\{i_1\}$, after which a corresponding point in the central subset $J_c\setminus\{j_1\}$ of the cloud $Y$ is chosen such that it lies close to the extreme circle and maximizes the residual component $|R^1_{j,i_2}|$. The algorithm refines the selection by iteratively searching for the best candidate, balancing proximity and error minimization, until the optimal point is found. 
As can be seen, this algorithm combines geometric and algebraic considerations and has a linear complexity. Moreover, contrary to the classical ACA it requires evaluation of only a small fraction of $R^1_{ij}$ components and thus could be even more advantageous than the classical ACA, especially for situations where the evaluation of the kernel is expensive. The complexity of the algorithm remains $\mathcal O(|Y|)$, but the number of operations is of the order $\beta |J^c| \approx \beta \lfloor\varepsilon_r |Y|-1\rfloor$ with the factor $\beta \sim 2\text{--}4$.

\begin{algorithm}
  \caption{Selection of optimal points for rank $k=2$}
  \label{alg:rank2_selection}
  \begin{algorithmic}[1]
    \Require Central subsets $I^c, J^c$ (excluding the first pivot)
      \State Select a point $x_{i_2}$ randomly such that $i_2 \in I^c$.
      \State Initialize $l \gets 1$.
      \State Copy the central subset $J^c_l \gets J^c$.
      \Repeat
          \State Select a point $y_{j^l_2}$ such that $j^l_2 \in J^c_l$ and minimizes the distance to the circle $\mathcal{C}_2(x_{i_2})$:
          \begin{equation*}
              j_2^l \gets \arg\min_{j_2 \in J^c_l} \mathrm{dist}(y_{j_2}, \mathcal{C}_2(x_{i_2}))
          \end{equation*}
          \State Evaluate $|R^1_{i_2,j_2^l}| = |a_{i_2,j_2^l} - a^{'}_{1(i_2,j_2^l)}|$.
          \If{$l>1$ \textbf{and} $|R^1_{i_2,j_2^l}| \le |R^1_{i_2,j_2^{l-1}}|$}
              \State Stop iterations and select $j_2 = j_2^{l-1}$.
          \Else
              \State $l \gets l + 1$
              \State Update the copy of the central subset $J^c_l \gets J^c_l \setminus \{j_2^l\}$ 
          \EndIf
      \Until convergence
  \end{algorithmic}
  \end{algorithm}

\subsubsection{Rank $k=3$}

\begin{figure}[ht!]
  \includegraphics[width=1\textwidth]{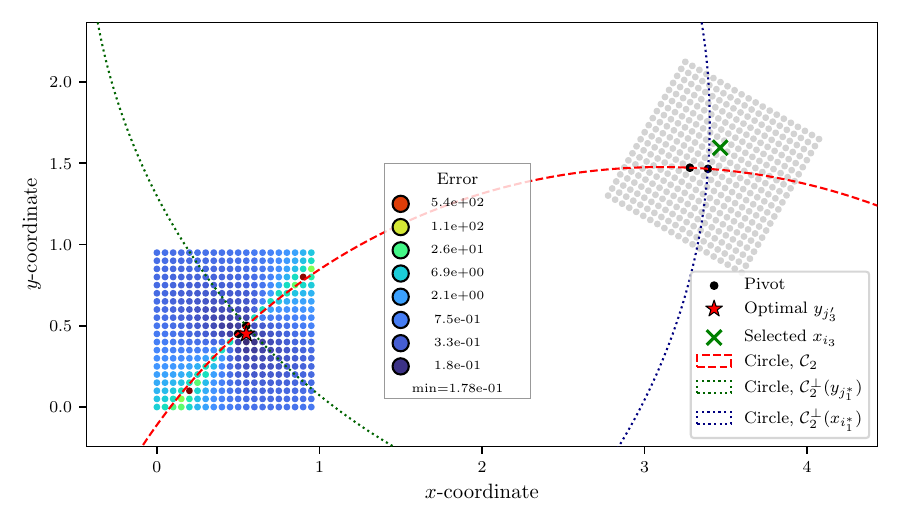}
  \includegraphics[width=1\textwidth]{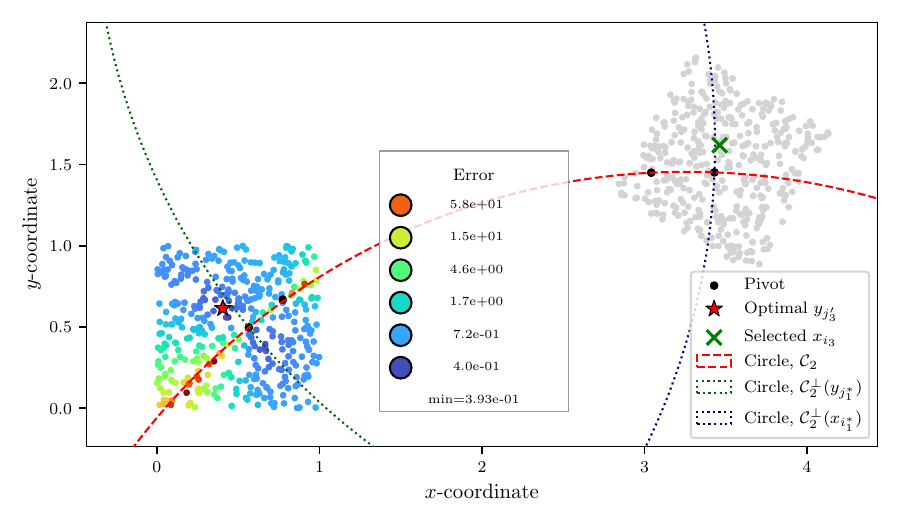}
  \caption{\label{fig:rank3good}Error structure for a random choice of $i_3$ for a structured (upper panel) and random (lower panel) clouds.}
\end{figure}

\begin{figure}[ht!]
  \includegraphics[width=1\textwidth]{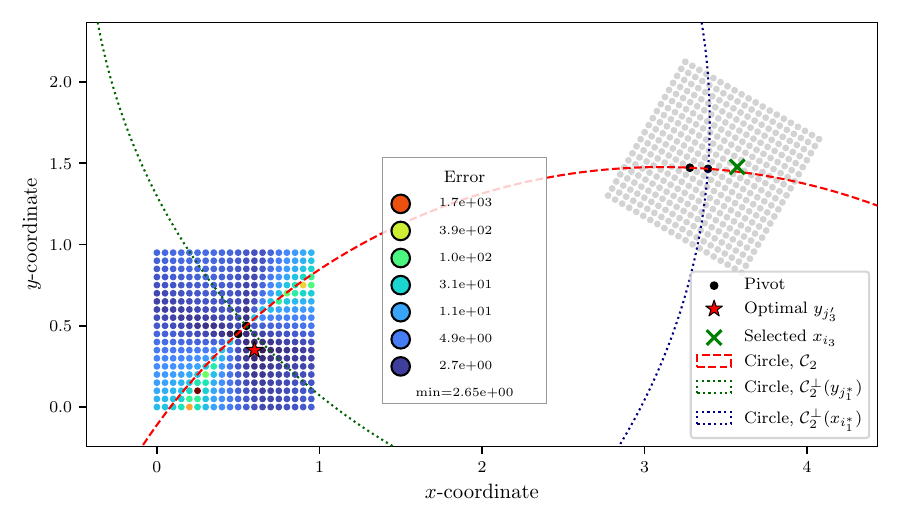}
  \includegraphics[width=1\textwidth]{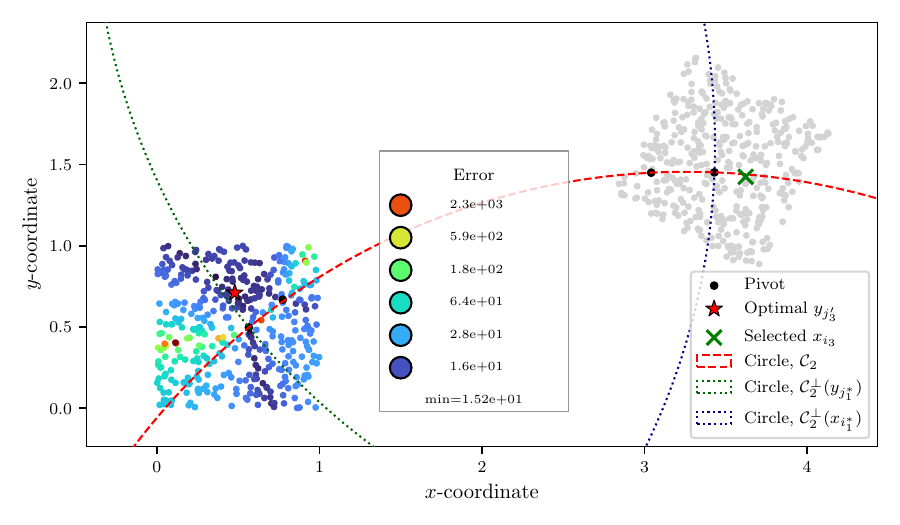}
  \caption{\label{fig:rank3bad}An example, of the error structure for a bad choice of $x_{i_3}$, i.e. near the extreme circle $\mathcal C_2(i^*_2)$ resulting in a significant increase in the error: upper panel for structured and lower panel for random clouds.}
\end{figure}

Third rank approximation is more trivial than the second one because the underlying geometry of errors remains stable with respect to the choice of the $x_{i_3}$ point under the condition that the latter is not located too close to $\mathcal C_2(i^*_2)$ which becomes now the extreme curve of maximal errors $\mathcal E_3^{\max}$. On the other hand, the conjugate extreme surface $\mathcal C^{\perp}_2(j_1^*)$ becomes the minimal error curve $\mathcal E_3^{\min}$, i.e. the extreme curves exchange their roles:
$$
\mathcal E_3^{\max} = \mathcal E_2^{\min} = \mathcal C_2(i^*_2), \quad \mathcal E_3^{\min} = \mathcal E_2^{\max} = \mathcal C_2^{\perp}(j^*_1).
$$
The error structure shown in \Cref{fig:rank3good} remains practically independent of the choice of $i_3$ if it is selected far from the $\mathcal E_3^{\max}$ extreme curve (see supplementary material for more examples, \cite{supplementary}). 
On the contrary, the selection of $x_{i_3}$ close to $\mathcal C_2(i^*_2)$ as shown in \Cref{fig:rank3bad} leads to a significant increase in the error (see the minimal value in the figure) and can alter the geometry of the error distribution.
Therefore, we assume that if one selects $x_{i_3}$ far from $\mathcal C_2(i^*_2)$, i.e. close to the circle passing through $x_{i^*_1}$ and orthogonal to the circle $\mathcal C_2(i^*_2)$: $\mathcal C_2^{\perp}(i^*_1)$ (our assumption), and then one selects $y_{j_3}$ close to $\mathcal C_2^{\perp}(j^*_1)$, then the error will be minimal.

The \Cref{alg:rank3_selection} for selecting the points $i_3$ and $j_3$ is rather similar to the one used for $i_2,j_2$. The main difference is that the $i_3$ point is selected in the central subset close to the $\mathcal C_2^{\perp}(x_{i^*_1})$ and $j_3$ is selected in an iterative way close to $\mathcal C_2^{\perp}(y_{j^*_1})$. The algorithm is as follows:
\begin{itemize}
  \item Select $i_3 \in I^c\setminus\{i_1,i_2\}$ such that 
  $$i_3 = \arg\min_{i\in I^c\setminus\{i_1,i_2\}} \mathrm{dist}(x_{i},\mathcal C_2^{\perp}(x_{i^*_1})).$$
  \item Select trial $j^1_3 \in J^c\setminus\{j_1,j_2\}$ such that
  $$j^1_3 = \arg\min_{j\in J^c\setminus\{j_1,j_2\}} \mathrm{dist}(y_{j},\mathcal C_2^{\perp}(y_{j^*_1})).$$
  \item Evaluate $|R^2_{i_3,j^1_3}| = |a_{i_3,j^1_3} - a'^{2}_{i_3,j^1_3}|$.
  \item Iterate to find $j_3$ such that
  $$j^{l+1}_3 = \arg\min_{j\in J^c\setminus\{j_1,j_2,j^1_3,\dots,j^l_3\}} \mathrm{dist}(y_{j},\mathcal C_2^{\perp}(j^*_1))$$
  before $|R^2_{i_3,j^{l+1}_3}| \le |R^2_{i_3,j^{l}_3}|$ and select $j_3 = j^{l}_3$ as the optimal point.
\end{itemize}
The complexity of the algorithm is $\mathcal O(|X|+|Y|)$, but the number of operations is of the order $|I^c|+\beta |J^c| \approx \lfloor\varepsilon_r |X|-2\rfloor + \beta \lfloor\varepsilon_r |Y|-2\rfloor$ with the factor $\beta \sim 2\text{--}4$.
Note that, contrary to the second rank selection, this algorithm for the selection of the 3rd rank pivot is not restricted to square-shaped domains and, in principle, can be applied to arbitrary domains.

\begin{algorithm}
  \caption{Rank $k=3$ Selection Algorithm}
  \label{alg:rank3_selection}
  \begin{algorithmic}[1]
    \Require Central subsets $I^c,J^c$ (excluding previously selected points)
    \State Select a point $i_3$ such that
    \begin{equation*}
      i_3 \gets \arg\min_{i \in I^c} \mathrm{dist}(x_i, \mathcal{C}_2^{\perp}(x_{i^*_1}))
    \end{equation*}
    \State Initialize $l \gets 1$
    \State Copy the central subset $J^c_l \gets J^c $
    \Repeat
      \State Select a point $j_3^l$ such that
      \begin{equation*}
        j_3^l \gets \arg\min_{j \in J^c_l} \mathrm{dist}(y_j, \mathcal{C}_2^{\perp}(y_{j^*_1}))
      \end{equation*}
      \State Evaluate 
      $|R^2_{i_3, j_3^l}| = |a_{i_3, j_3^l} - a'_{2(i_3,j_3^l)}|.$
      \If{$l>1$ \textbf{and} $|R^2_{i_3, j_3^l}| \le |R^2_{i_3, j_3^{l-1}}|$}
        \State Stop iterations and select $j_3 = j_3^{l-1}$
      \Else
        \State $l \gets l + 1$
        \State Update the copy of the central subset $J^c_l \gets J^c_l \setminus \{j_3^l\}$
      \EndIf
    \Until convergence
  \end{algorithmic}
\end{algorithm}

\subsubsection{Higher ranks $k\ge 4$\label{sec:algorithm_for_high_ranks}}

For higher ranks establishing a simple deterministic procedure relying on extreme curves seems too tedious: trivial curves do not ensure extreme values of error even though visually it can be noticed that they also form curves with more or less constant curvature. Therefore, an even simpler approach combining geometrical consideration and algebraic analysis could be used. The algorithm is presented in 
\Cref{alg:higher_ranks} and consists in a two step evaluation of the maximal residual component over the central subsets. As previously, we introduce adjusted central subsets by excluding previosly selected points:
$$ I^c_k = I^c \setminus \{i_1,\dots,i_{k-1}\}, \quad J^c_k = J^c \setminus \{j_1,\dots,j_{k-1}\}.$$
Then for a randomly selected trial $i^t \in I^c_k$ we evaluate the components of the residual matrix corresponding to the columns $j \in J^c_k$:
$$
R^{k-1}_{i^tj} = a_{i^tj} - \sum_{l=1}^{k-1} u_l(i^t) v_l(j).
$$
Among all evaluated components we select the column $j_k$ with the maximal absolute value of the residual component:
$$
j_k = \arg\max_{j \in J^c_k} |R^{k-1}_{i^tj}|.
$$
Then we evaluate the components of the central rows $i \in I^c_k$:
$$
R^{k-1}_{ij_k} = a_{ij_k} - \sum_{l=1}^{k-1} u_l(i) v_l(j_k).
$$
and select the row $i_k$ with the maximal absolute value of the residual component:
$$
i_k = \arg\max_{i \in I^c_k} |R^{k-1}_{ij_k}|.
$$
In such a way, the new pivot $i_k,j_k$ is selected. This algorithm will be used for all ranks $k\ge 4$.  

\begin{algorithm}
\caption{Pivot Choice for Higher Ranks $k \ge 4$}
\label{alg:higher_ranks}
\begin{algorithmic}[1]
  \Require Central subsets $I^c,J^c$ (excluding previously selected points)
  \State Select a trial point $i^t$ randomly such that $i^t \in I^c_k$
  \State Compute $R^{k-1}_{i^t j}$ for all $j \in J^c$:
  \begin{equation*}
    R^{k-1}_{i^t j} = a_{i^t j} - \sum_{l=1}^{k-1} u_l(i^t) v_l(j)
  \end{equation*}
  \State Select the column $j_k$ such that
  \begin{equation*}
    j_k \gets \arg\max_{j \in J^c} |R^{k-1}_{i^t j}|
  \end{equation*}
  \State Compute $R^{k-1}_{i j_k}$ for all $i \in I^c$:
  \begin{equation*}
    R^{k-1}_{i j_k} = a_{i j_k} - \sum_{l=1}^{k-1} u_l(i) v_l(j_k)
  \end{equation*}
  \State Select the row $i_k$ with the maximal absolute value:
  \begin{equation*}
    i_k \gets \arg\max_{i \in I^c_k} |R^{k-1}_{i j_k}|
  \end{equation*}
\end{algorithmic}
\end{algorithm}

\subsection{A New Method: ACA-GP}

The previous sections suggest a methodology for an informed choice of pivots which aim at minimizing the approximation errors. 
In this chapter the full algorithm for a new adaptive low-rank approximation is presented and discussed.
In essence, we suggest combining the purely algebraic ACA method with geometrical considerations for the choice of the pivots. This choice is based on combined knowledge of spatial distribution of points within the cloud and on the algebraic properties of the residual vectors. The proposed method is named ACA-GP (GP for Geometrical Pivots). 

\begin{algorithm}[ht!]
  \caption{Adaptive Cross Approximation with Geometrical Pivots (ACA-GP) for $m\le n$\label{alg:acagp}}
  \begin{algorithmic}[1]
  \Require Two clouds $X = \{x_i\}_{i=1}^n$, $Y = \{y_j\}_{j=1}^m$, Kernel $\kappa(x,y)$, global tolerance $\epsilon >0$ \textbf{OR} maximal rank $k_{\max}$ \newline pivot tolerance $\epsilon_p > 0$, central fraction $\varepsilon_r > 0$
  \Statex \textbf{// First pivot}
  \State Find the center of $X$: $\bar x \gets 1/n\;\sum_{i\in N} x_i$ \Comment $\mathcal O(n)$
  \State Find the center of $Y$: $\bar y \gets 1/m\;\sum_{j\in M} y_j$ \Comment $\mathcal O(m)$
  \State Find the first pivot's row: $\displaystyle i_1 \gets \arg\min_{\substack{i \in N \\ (x_i-\bar x)\cdot (\bar y-\bar x) > 0}} \|x_i - \bar{x}\|$  \Comment $\mathcal O(n)$\label{alg:acagp:i1}
  \State Find the first pivot's column: $\displaystyle j_1 \gets \arg\min_{\substack{j \in M \\(y_j- \bar y)\cdot (\bar x - \bar y) > 0}} \|y_j - \bar{y}\|$ \Comment $\mathcal O(m)$\label{alg:acagp:j1}
  \State $I \gets \{i_1\}$, $J \gets \{j_1\}$.
  \State Pivot $p_1 = \kappa(x_{i_1}, y_{j_1})$
  \State $u^1 \gets \{\sign(p_1)\kappa(x_{i}, y_{j_1}) / \sqrt{p_1}\}, i\in N$,  \Comment $\mathcal O(m)$\label{alg:acagp:step_u1}
  \State $v^1 \gets \{\kappa(x_{i_1}, y_{j}) / \sqrt{p_1}\}, i \in M$, \Comment $\mathcal O(n)$\label{alg:acagp:step_v1}
  \State $U \gets u^1$, $V \gets v^1$
  \State Compute residuals: $\|R\|_F = \|U\|_F\|V\|_F$, $\|A'\|_F = \|R\|_F$ \Comment $\mathcal O(n+m)$
  \State Construct central subsets $I_c,J_c$ according to \Cref{alg:central_subset}.
  \Statex \textbf{// Main loop}
  \State $r \gets 2$
  \Repeat
      \If {$r = 2$}  \State Use \Cref{alg:rank2_selection} \textbf{or} \Cref{alg:higher_ranks} to select $i_2,j_2$ \Comment $\mathcal O(m)$ or $\mathcal O(n+m)$
      \ElsIf {$r = 3$} \State Use \Cref{alg:rank3_selection} to select $i_3,j_3$ \Comment $\mathcal O(m)$
      \Else \State Use \Cref{alg:higher_ranks} to select $i_r,j_r$ \Comment $\mathcal O(n+m)$
      \EndIf
      \State Update central subsets $I_c \gets I_c\setminus\{i_r\}$, $J_c \gets J_c\setminus\{j_r\}$
      \State Update pivot sets $I \gets I\cup\{i_r\}$, $J \gets J\cup\{j_r\}$
      \State Pivot $p_r = \kappa(x_{i_r}, y_{j_r}) - \sum_{l=1}^{r-1} u_l(i_r) v_l(j_r)$
      \If {$|p_r| < \epsilon_p$}
      \Return $U,V,k \gets r-1$ 
      \EndIf
      \State Evaluate column $\tilde{u}^r := \kappa(x_i, y_{j_r}), i\in N$ \Comment $\mathcal O(n)$
      \State Evaluate row $\tilde{v}^r := \kappa(x_{i_r}, y_j), j\in M$ \Comment $\mathcal O(m)$
      \For{$l = 1, \ldots, r - 1$} \Comment Substract $A'_{r-1}$ columns and rows, $\mathcal O(rn)$\label{alg:acagp:step_x}
        \State $\tilde u^r := \tilde u^r - (v_l)_{j_r} u_l$ 
        \State $\tilde v^r := \tilde v^r - (u_l)_{i_r} v_l$
      \EndFor
      \State Compute skeleton row and column: $u^r \gets \sign(p_r)\tilde{u}^r / \sqrt{p_r}$, $v^r \gets \tilde{v}^r / \sqrt{p_r}$
      \State Update outer-product matrices: $U \gets [U|u^r]$, $V \gets [V|v^r]$
      \State Residual norm $\|R_r\|_F = \|u^r\|_F \|v^r\|_F$ \Comment $\mathcal O(n+m)$
      \State Compute the norm of the approx. matrix $\|A'_r\|_F$ \Cref{eq:mnorm} \Comment $\mathcal O(r(n+m))$\label{alg:acagp:mnorm}
      \State $r \gets r + 1$
  \Until{$\|R_r\|_F \leq \epsilon \|A'_r\|_F$ or $J \equiv M$ or $r = k_{\max}$}
  \end{algorithmic}
\end{algorithm}

A detailed algorithm is presented in~\Cref{alg:acagp} with associated complexity evaluations for every step. We assumed that $m\le n$. The following notations were used: two sets of integers $M=\{1,2,\dots,m\}$, 
$N = \{1,2,\dots,n\}$ indicate IDs of points (elements) in two separate clouds, whereas sets $I_k=\{i_1,i_2,\dots,i_k\}$, $J_k = \{j_1,j_2,\dots,j_k\}$
 correspond to pivots (rows and columns) used for skeleton construction of a rank-$k$ approximation, 
so $N\setminus I_k$, $M\setminus J_k$ correspond to the sets of the remaining columns and rows, respectively. 
The sets $I^c,J^c$ are the central subsets which are defined by \Cref{eq:central_subsets}.
The core of the method remains unchanged with respect to the ACA methodology, it is briefly described below to help the interpretation of the algorithm.
The matrix $A'_{k}$ mentioned in \Cref{alg:acagp:step_x} corresponds to the matrix approximation at the $k$-th step/rank and is defined as $A'_k = U_k V_k^T$ with $U_k = \{u^1,u^2,\dots,u^k\}$ and $V_k = \{v^1,v^2,\dots,v^k\}$, but of course $A'_k$ is never evaluated, only matrices $U_k,V_k$ are fully computed and only the necessary operations are performed to construct the next approximation $A'_{k+1}$. The rows and columns $u^i,v^i$ for $i = 1,\dots,k$ of the matrices $U_k,V_k$ are constructed from scaled rows and columns of the residual matrices $R^{i-1}$ with $R^0 = A$ in the following way. For the first rank $k=1$ the rows and columns are constructed as
$$u^1_i = \sign(p_1) A_{i,j_1}/\sqrt{p_1}, \quad v^1_j = A_{i_1,j} /\sqrt{p_1},$$ 
where $p_1 = A_{i_1,j_1}$ is the first pivot (see \Cref{alg:acagp:step_u1,alg:acagp:step_v1}). The division by pivot is split into two divisions for $u$ and $v$ vectors to avoid division by a too small number.
For higher ranks $k>1$, the rows and columns are constructed as
$$
u^k_i = \sign(p_k) \tilde{u}^k_i/\sqrt{p_k}, \quad v^k_j = \tilde{v}^k_j/\sqrt{p_k},
$$
where
$$
\tilde u^k_i = R^{k-1}_{ij_k}  = a_{i,j_k} - \sum_{l=1}^{k-1} v^l_{j_k} u^l_i, \quad \tilde v^k_j = R^{k-1}_{i_k,j} = a_{i_k,j} - \sum_{l=1}^{k-1} u^l_{i_k} v^l_j,
$$
where $p^k = \tilde{u}^k_{i_k} = \tilde{v}^k_{j_k}$ is the $k$-th pivot. 

To evaluate the Frobenius norm of the approximate matrix at the $k$-th step/rank (lines 35), we use the recursive formula from~\parencite{bebendorf2003adaptive}:
\begin{equation}
  \|A'_{k}\|_F = \sqrt{\|A'_{k-1}\|^2_F + 2 \sum\limits_{j=1}^{k-1} u_k^\top u_j v_j^\top v_k + R^2_k}, \quad \text{ where } \quad R^2_k = \|u_k\|^2_F \|v_k\|^2_F 
  \label{eq:mnorm}
\end{equation}

The key steps of the entire procedure is the selection of to-be-evaluated rows and columns of the original matrix $A$, which is done according to algorithms presented in the previous section, namely \Cref{alg:rank2_selection,alg:rank3_selection,alg:higher_ranks}.
The selection of the first pivot is made such that the pivot is the closest point to the barycenter of the cloud and lies in the half-plane passing through the barycenter on the side containing the homologue cloud, see \Cref{alg:acagp:i1,alg:acagp:j1}:
\begin{align}
  i_1 &= \arg \min_{i\in N,(x_i-\bar x)\cdot (\bar y-\bar x) \ge 0} \|x_i - \bar x\|\\ 
  j_1 &= \arg \min_{j\in M,(y_j-\bar y)\cdot (\bar x-\bar y) \ge 0} \|y_j - \bar y\|,
\end{align}
where $\bar x, \bar y$ are the barycenters of clouds $X,Y$, respectively.

\section{Results: ACA-GP for Two Interacting Clouds\label{sec:results}}

\begin{figure}[ht!]
  \centering\includegraphics[width=0.75\textwidth]{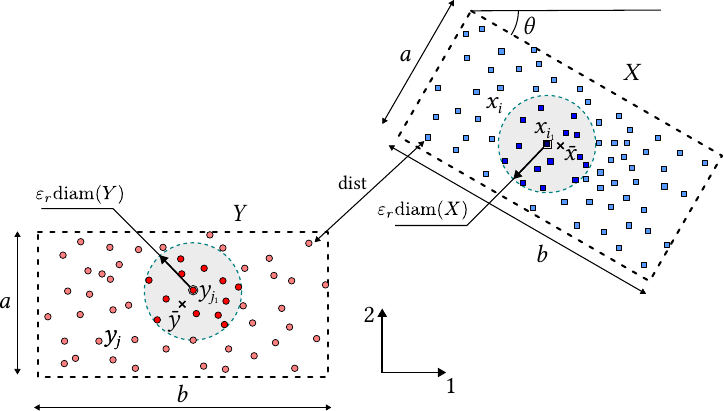}
  \caption{\label{fig:exp_setup}Set-up of the experiment: two clouds of points $X$ and $Y$ are generated, each in a rectangle $a\times b$ ($\xi = a/b$).}
\end{figure}

We demonstrate the performance of the algorithm on the problem of two interacting separate clouds of points $X = \{x_i\} \subset \mathbb{R}^{n \times d}$, $Y = \{y_j\} \subset \mathbb{R}^{m \times d}$ where $d=2$ is the space dimension. 
For the sake of simplicity, points in every cloud $X,Y$ are generated 
within a rectangle $a\times b$ such that $a\le b$, and $\xi = a/b$ is the aspect ratio of the rectangle, $\xi \in (0,1]$ (see Fig.~\ref{fig:exp_setup}). 
Therefore, the size of every cloud can be approximated by $\mathrm{diam}(Y) \approx \sqrt{a^2+b^2} = b\sqrt{1+\xi^2}$.
For simplicity, we set $b=1$, so the distances are measured in units of $b$.
We consider a Poisson point process (uniform distribution) for both clouds $X,Y$ to generate the points.
The cloud $Y$ is centred at zero and the rectangle is aligned with the axes, whereas the cloud $X$ is displaced and rotated by a random angle $\theta$ (uniform distribution in $[-\pi,\pi]$) in such a way that the criterion of admissibility \Cref{eq:admissibility} is satisfied. The approximate distance between the clouds is estimated as 
\begin{equation}
  \mathrm{dist}'(X,Y) = \|\bar x - \bar y\| - \min\{\mathrm{diam}(X), \mathrm{diam}(Y)\} = \|\bar x - \bar y\| - b\sqrt{1+\xi^2},
\end{equation}
whereas the true distance can be computed as
\begin{equation}
  \mathrm{dist}(X,Y) = \min\limits_{i \in N,j \in M} \|x_i - y_j\|.
\end{equation}
The kernel is given by a smooth kernel $\kappa(x,y) = 1/\|x-y\|$ and we assume that the matrix entries are simply
$$
a_{ij} = \kappa(x_i, y_j)
$$
as it was considered in~\parencite{grasedyck2005adaptive}.

The tests are carried out for two aspect ratio $\xi = \{0.5, 1\}$ and for true distances $\mathrm{dist} = \{1.5, 2.5, 5.0\}$, the number of points in the clouds is $n,m = 400$.
We used $N_s=1000$ realization of randomly oriented clouds with randomly distributed points. The first objective is to study the effect of the central subset fraction $\varepsilon_r$ introduced in \Cref{sec:central_subsets}.
The rule of thumb for its selection can be easily constructed. 
If we suppose that the domain is a circle of diameter $\text{diam}(Y)$, so to have the number of points in the central subset of radius $\varepsilon_r \text{diam}(Y)$ equal to $\max(\text{rank})$, we need to verify the condition:
$$
\varepsilon_r \gtrsim \sqrt{\frac{\max(\text{rank})}{|Y|}},
$$
which is a double of the non-conservative estimation for the needed value. The double is used to account for arbitrary shape of the domain and the fact that the number of points in the central subset should be bigger than the maximal rank of the approximation to ensure optimal performance. So, one has to carefully choose the maximal allowed rank in the algorithm to avoid a too large central subset fraction and excessive computational cost. In the current case for $\max(\text{rank})=10$, we get $\varepsilon_r \gtrsim 0.15$.

\begin{figure}[ht!]
  \includegraphics[width=1\textwidth]{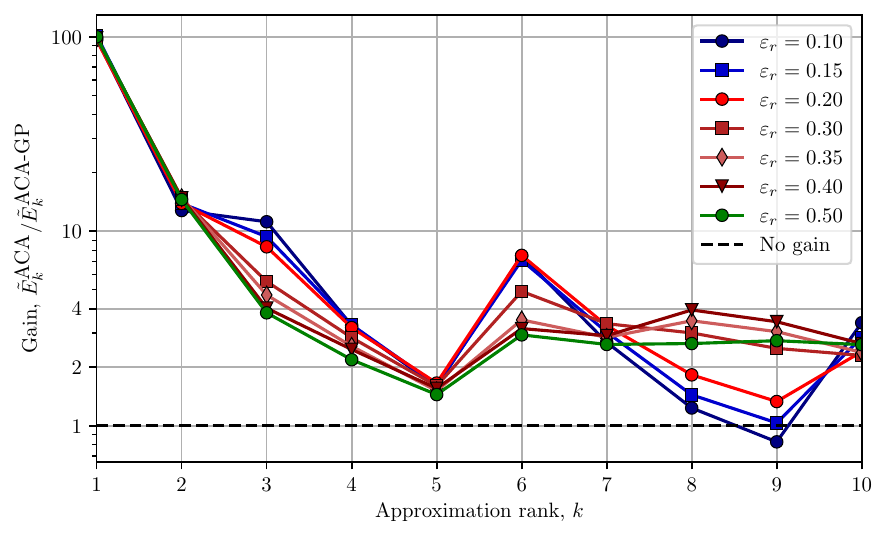}
  \caption{\label{fig:effect_of_center}The effect of the relative central subset fraction $\varepsilon_r$ on the accuracy of the low-rank approximation: the gain factor $\tilde E_k^{\textrm{ACA-GP}}/\tilde E_k^{\textrm{ACA}}$ for all ranks are plotted.}
\end{figure}

The effect of the central subset fraction $\varepsilon_r$ is not trivial and affects differently different ranks of the approximation. 
The relative size of the central subset is varied in the interval $\varepsilon_r \in [0.1, 0.5]$. 
The lower value was chosen such to ensure that a sufficient number of points is available in each central subset even though in Poisson's process it is always possible to have zero points within a central subset. 
In \Cref{fig:effect_of_center} we show for $\xi=1$ and $\mathrm{dist}=1.5$ how the accuracy of the approximation varies with the central subset size for the first 10 ranks.
The very first rank is, of course, independent on the central ratio, the second rank is also neutral to the variation of the parameter $\varepsilon_r$. 
Starting from rank $k=3$ and up to $k=7$, the increase of $\varepsilon_r$ decreases the gain factor, but for ranks $k=8,9$ the gain factor is higher for larger $\varepsilon_r$.
However, a simple strategy consisting in increasing the central subset individually for those ranks does not work because the overall approximation depends on the whole history of pivot selections. Therefore, a selected central subset fraction should be kept for all ranks to ensure a controllable effect.
The value of $\varepsilon_r=0.25$ presents a reasonable tradeoff between the gain factor for all considered ranks. However, if one is interested in the lower ranks only (up to $k=6$), one can choose the smallest $\varepsilon_r=0.1$ to ensure the best accuracy and less computational cost.

\begin{figure}
  \begin{subfigure}{0.48\textwidth}
    \includegraphics[width=1\textwidth]{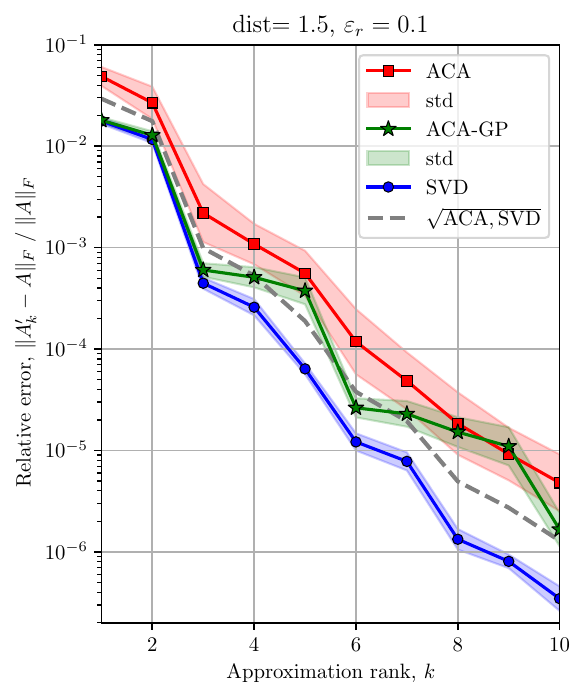}
    \caption{\label{fig:aca_gp_error_xi_1_dist_15}}
  \end{subfigure}
  \begin{subfigure}{0.48\textwidth}
    \includegraphics[width=1\textwidth]{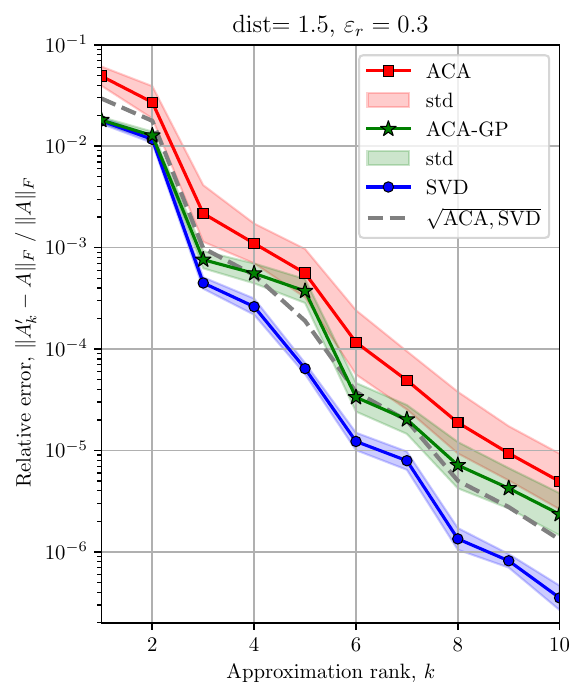}
    \caption{\label{fig:aca_gp_gain_xi_1_dist_15}}
  \end{subfigure}
  \begin{subfigure}{0.48\textwidth}
    \includegraphics[width=1\textwidth]{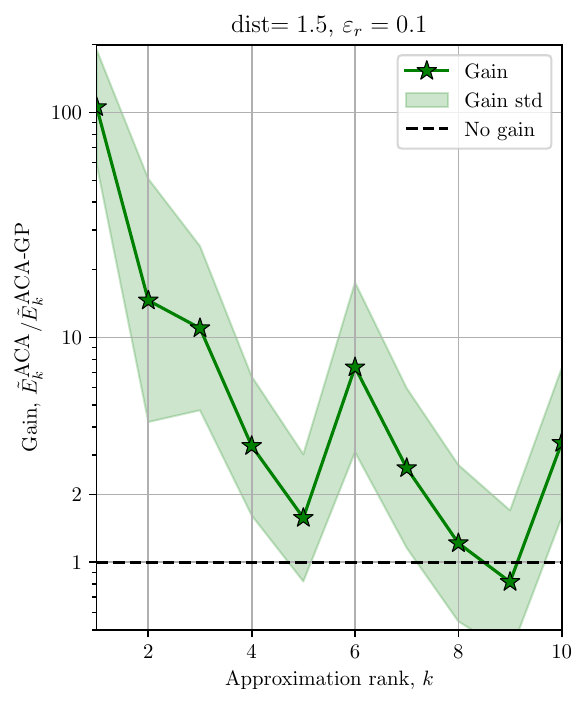}
    \caption{\label{fig:aca_gp_error_xi_1_dist_25}}
  \end{subfigure}
  \begin{subfigure}{0.48\textwidth}
    \includegraphics[width=1\textwidth]{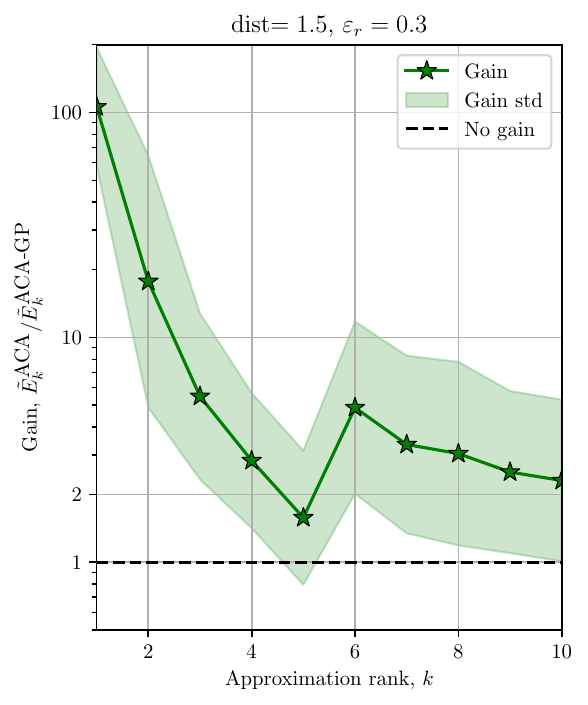}
    \caption{\label{fig:aca_gp_gain_xi_1_dist_25}}
  \end{subfigure}
  \caption{\label{fig:aca_gp_error_xi_1}(a-b) relative error in matrix approximation by ACA, ACA-GP and SVD methods for square domains ($\xi=1$): log-mean values and log-standard deviations are computed over 1000 realizations; (c-d) corresponding log-mean gain and log-standard deviation between the ACA-GP compared to the ACA method with the offset of SVD accuracy (\Cref{eq:gain}) is also shown; (a,c) corresponds to $\varepsilon_r=0.1$ and (b,d) to $\varepsilon_r=0.3$. The true distance is $\mathrm{dist}(X,Y)=1.5$.}
\end{figure}

\begin{figure}
  \begin{subfigure}{0.48\textwidth}
    \includegraphics[width=1\textwidth]{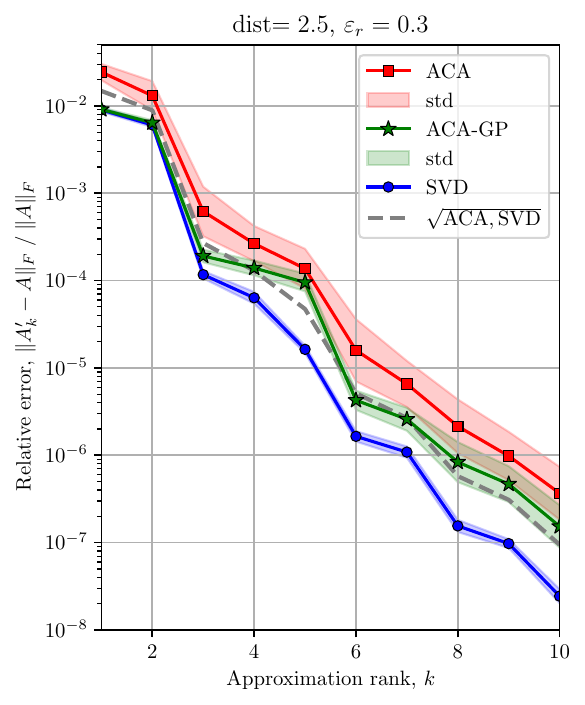}
    \caption{\label{fig:aca_gp_error_xi_05_dist_15}}
  \end{subfigure}
  \begin{subfigure}{0.48\textwidth}
    \includegraphics[width=1\textwidth]{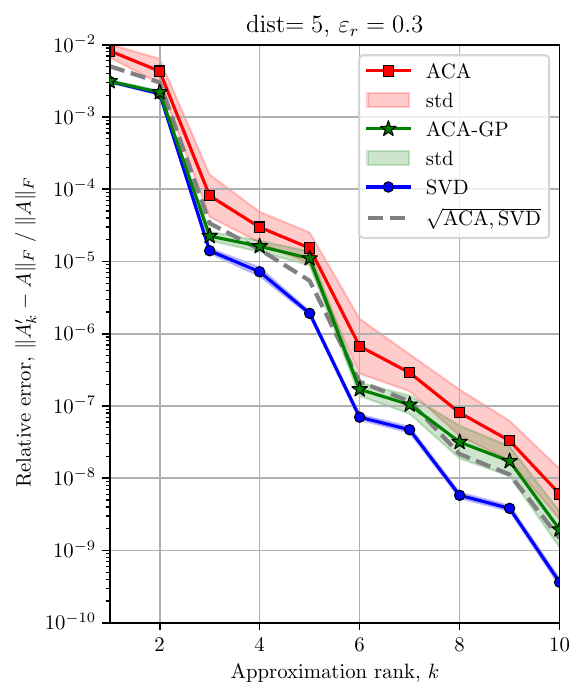}
    \caption{\label{fig:aca_gp_gain_xi_05_dist_15}}
  \end{subfigure}
  \begin{subfigure}{0.48\textwidth}
    \includegraphics[width=1\textwidth]{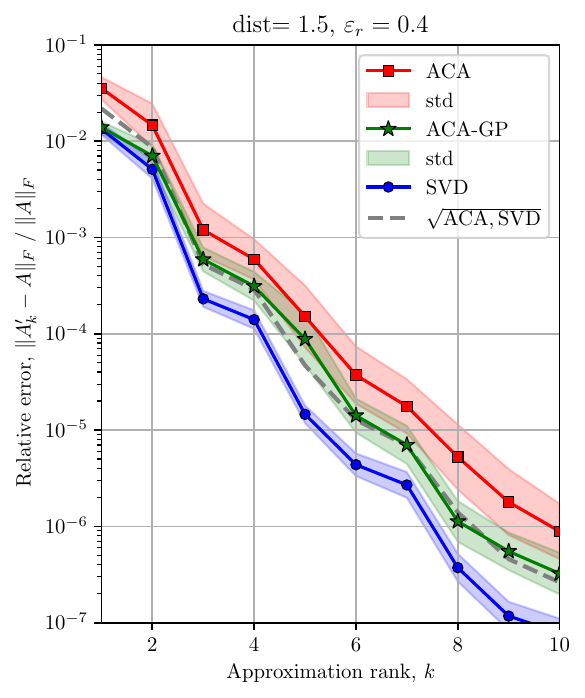}
    \caption{\label{fig:aca_gp_error_xi_05_dist_25}}
  \end{subfigure}
  \begin{subfigure}{0.48\textwidth}
    \includegraphics[width=1\textwidth]{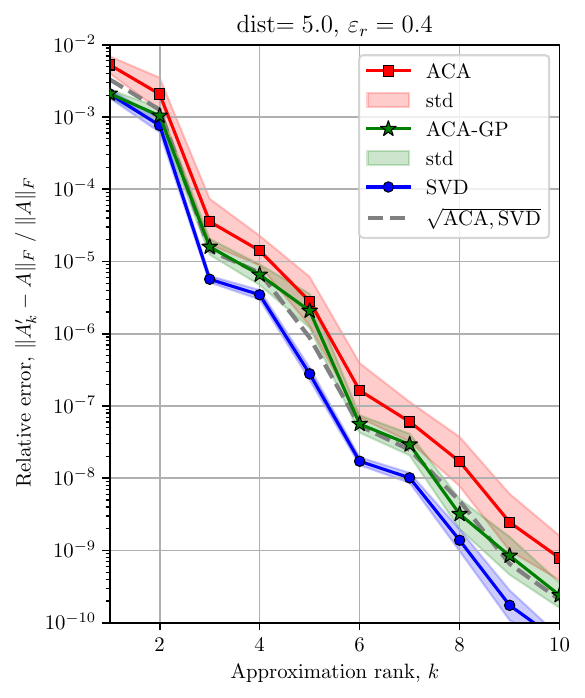}
    \caption{\label{fig:aca_gp_gain_xi_05_dist_25}}
  \end{subfigure}
  \caption{\label{fig:aca_gp_error_xi_05}Relative error in matrix approximation by ACA, ACA-GP and SVD methods for square domains $\xi=1$ (computed over 1000 realizations) for $\mathrm{dist}=2.5$ in (a) and $\mathrm{dist} = 5.0$ in (b) for $\varepsilon=0.3$;
  relative errors for rectangular domains $\xi=0.5$ (computed over 500 realizations) are shown in (c-d) for $\mathrm{dist}=1.5$ and $\mathrm{dist}=5.0$ with $\varepsilon=0.4$.}
\end{figure}

The main results are presented in \Cref{fig:aca_gp_error_xi_1,fig:aca_gp_error_xi_05} for different aspect ratios $\xi$, true distances $\mathrm{dist}$ and central subset sizes $\varepsilon_r$. We plot the relative Frobenius norm of the residual matrix $\|A-A'\|_F/\|A\|_F$ for the three methods: ACA, ACA-GP and SVD; the log-mean values and log-standard deviations are identified as:
\begin{align}
E_{\log} &= \log_{10} \frac{\|A - A'_k\|_F}{\|A\|_F}, \quad \bar E_{\log} = \frac{1}{N_s} \sum_{s=1}^{N_s} E_{\log}^{(s)}\\
\sigma_{\log} &= \sqrt{\frac{1}{N_s} \sum_{s=1}^{N_s} (E_{\log}^{(s)} - \bar E_{\log})^2}, \quad E = 10^{\bar E_{\log} \pm \sigma_{\log}}.
\end{align}
The accuracy gain factor of the ACA-GP to the ACA compared to the SVD approximation is also plotted:
\begin{equation}
  \text{gain} = \frac{\|A - A'_{\text{ACA}}\|_F - \|A - A_{\text{SVD}}\|_F}{\|A - A'_{\text{ACA-GP}}\|_F - \|A - A_{\text{SVD}}\|_F}.
  \label{eq:gain}
\end{equation}
For aspect ratio $\xi=1$, the results are computed over 1000 realizations, for $\xi=0.5$ the results are computed over 500 realizations. Many more results can be found in the supplementary material~\parencite{supplementary}.
The key observations are the following:
\begin{itemize}
  \item The ACA-GP method strongly outperforms the classical ACA method for most of the considered ranks.
  \item The average error of the ACA-GP method and its standard deviation are significantly lower than those of the ACA method.
  \item For square-shaped domains and small $\varepsilon_r$, the first 3 ranks, the ACA-GP method matches very closely in its performance with the SVD approximation if the extreme curves are used through \Cref{alg:rank2_selection,alg:rank3_selection}.
  \item For higher ranks $k\ge4$ and arbitrary domains (here, rectangular) and properly selected $\varepsilon_r$, the ACA-GP's error represents the geometrical mean of the ACA and SVD errors.
  \item The ACA-GP is the most efficient up to rank $k=3$ and for rank $k=6$, for intermediate ranks $k=\{4,5\}$, the error decay rate remains small.
  \item The domain aspect ratio $\xi$ does not affect the performance of the ACA-GP method if the central subset-based selection is used \Cref{alg:higher_ranks}, eventually in combination with \Cref{alg:rank3_selection}.
\end{itemize}

\section{Conclusions\label{sec:conclusions}}

In the context of matrix construction for Boundary Integral or similar methods, the low-rank approximation of the interaction matrix is a key ingredient for the effective use of hierarchical matrices. Within this approach, the interaction between well separated or admissible domains can be efficiently approximated by low-rank matrices.
In this work, we proposed a new adaptive low-rank approximation method, called ACA-GP, which combines the algebraic ACA method with a geometry-aided choice of pivots. This new method, similarly to the ACA method, does not require the full original matrix to construct its low-rank approximation and relies on the algorithms with much lower complexity than the SVD. The new method demonstrates a much higher accuracy and has a smaller error dispersion than the classical ACA method.

The pivot selection procedure is the following.
For the first pivot, the central points/elements from the interacting domains are selected. The choice of the pivots for ranks 2 and 3 is based on specific  geometrical error structure revealed here, however, it can be well constructed for square-shaped domains only. Nevertheless, this structure is rather independent on the distribution of the points within the interacting domains.For higher ranks, adopting a simple geometric procedure, which constraints the choice of pivots to the central part of the domain, allows to improve the accuracy of the approximation and even reduce the computational cost compared to the classical ACA method. The method was tested on two randomly oriented interacting clouds of points of different aspect ratios. The performance of the ACA-GP method matches closely the SVD performance for the first three ranks (especially for square-shaped domains) and ensures approximation error equal in average to the geometric mean of the ACA and SVD for higher ranks.

In terms of the implementation, the ACA-GP method is slightly more complex and contains one user-defined parameter: the central subset fraction $\varepsilon_r$, which controls the size of the subdomain from which geometrical points corresponding to rows and columns should be selected. The optimal value of this parameter depends on the required rank of the approximation, which in its turn cannot be known in advance because of the key idea of the ACA -- the adaptivity -- the approximate matrix construction is stopped when the approximate error reaches a certain threshold. Nevertheless, keeping the central fraction as small as possible is beneficial for its performance for lower ranks, which should be often enough in practice. 
In the worst case scenario, the central subset can be increased when it runs out of available points.
Another possible drawback of the focus on the central subset is that in real life, the domains can have a ring-like structure, which therefore does not contain the central subset and could be tricky to handle. In such situation, it would be relevant to switch to the classical ACA method.

In perspective, we plan to adapt the ACA-GP for the collocation BIM method, where the area associated with each domain should be taken into account.
We also plan to study in detail the geometrical structure of extreme curves and make an attempt to suggest a simple algorithm for its identification for arbitrary domains.

Our Python implementation of the ACA and ACA-GP is shared on GitHub~\parencite{yastrebovGitHub} as well as the data and the scripts used for the experiments~\parencite{supplementary}.

\section*{Acknowledgement}

This work was partially presented at the World Congress of Computational Mechanics held in Vancouver in the summer of 2024, with fees partially covered by the Comité National Français de Mécanique (CNFM), which is gratefully acknowledged.




\end{document}